\documentclass [12pt,a4paper]{amsart} 

\usepackage[T1]{fontenc}

\usepackage{graphicx}
\usepackage[utf8]{inputenc} 
\usepackage{amssymb,amsmath,amsfonts,amsthm,mathtools}
\usepackage{enumitem}
\renewcommand{\raggedright}{\justifying}							
\usepackage{xcolor}
\usepackage{helvet}							
\usepackage{comment}

\usepackage{pgfplots}
\pgfplotsset{compat=1.15}
\usepackage{mathrsfs}
\usetikzlibrary{arrows}
\usepackage{relsize}

\newtheorem{thm}{Theorem}
\newtheorem{theorem}{Theorem}[section]
\newtheorem{lemma}[theorem]{Lemma}
\newtheorem{proposition}[theorem]{Proposition}
\newtheorem{corollary}[theorem]{Corollary}
\newtheorem{problem}{Problem}

\theoremstyle{definition}
\newtheorem{definition}[theorem]{Definition}

\theoremstyle{remark}

\newtheorem*{notation*}{Notation}
				
\newcommand{\hilbert}[3]{\left(\frac{{#1}, {#2}}{#3}\right)}
\newcommand{\Hyp}{\mathbb{H}}
\newcommand{\bH}{\mathbb{H}}
\newcommand{\cH}{\mathcal{H}}
\newcommand{\rM}{\mathrm{M}}

\newcommand{\rn}{\mathrm{n}}
\newcommand{\cO}{\mathcal{O}}
\newcommand{\rP}{\mathrm{P}}
\newcommand{\Q}{\mathbb{Q}}
\newcommand{\R}{\mathbb{R}}
\newcommand{\id}{\mathrm{id}}
\newcommand{\SL}{\mathrm{SL}}
\newcommand{\PSL}{\mathrm{PSL}}
\newcommand{\sys}{\mathrm{sys}}

\newcommand{\tr}{\mathrm{tr}}

\newcommand{\Mah}{\mathrm{M}}

\def\tareesidedbox#1{\setbox0=\hbox{$#1$}\dimen0=\wd0 \advance\dimen0 by3pt\rlap{\hbox{\vrule height9pt width.4pt depth2pt \kern-.4pt\vrule height9.4pt width\dimen0 depth-9pt\kern-.4pt \vrule height9pt width.4pt depth2pt}} \relax \hbox to\dimen0{\hss$#1$\hss}}
\def\ho#1{\tareesidedbox{#1}}

\usepackage{hyperref}

\usepackage{enumitem}
\title[Geometry and arithmetic of semi-arithmetic groups]{Geometry and arithmetic of semi-arithmetic Fuchsian groups}

\author[M. Belolipetsky]{Mikhail Belolipetsky}
\address{IMPA, Estrada Dona Castorina 110, 22460-320 Rio de Janeiro, Brazil}
\email[]{mbel@impa.br}

\author[G. Cosac]{Gregory Cosac} 
\address{Departamento de Matem\'atica Aplicada, IME - Universidade de S\~ao Paulo, Rua do Mat\~ao 1010, 05508-090 – São Paulo, SP, Brazil}
\email{cosac@ime.usp.br}

\author[C. D\'oria]{Cayo D\'oria}
\address{UFS, Departamento de Matem\'atica - Av. Marcelo D\'eda Chagas s/n, 49100-000. S\~ao Crist\'ov\~ao, Brazil}
\email{cayo@mat.ufs.br}

\author[G. Teixeira Paula]{Gisele Teixeira Paula}
\address{UFPR, Centro Politécnico - Av. Cel. Francisco H. dos Santos 100, 81530-000. Curitiba, Brazil}
\email[]{giseleteixeira@ufpr.br}

\begin{document}

\begin{abstract}
Semi-arithmetic Fuchsian groups is a wide class of discrete groups of isometries of the hyperbolic plane which includes arithmetic Fuchsian groups, hyperbolic triangle groups, groups admitting a modular embedding, and others. We introduce a new geometric invariant of a semi-arithmetic group called stretch. Its definition is based on the notion of the Riemannian center of mass developed by Karcher and collaborators. We show that there exist only finitely many conjugacy classes of semi-arithmetic groups with bounded arithmetic dimension, stretch and coarea. The proof of this result uses the arithmetic Margulis lemma. We also show that when stretch is not bounded there exist infinite sequences of such groups.
\end{abstract}



\maketitle 

\section{Introduction}

Semi-arithmetic Fuchsian groups were introduced by Schaller and Wolfart \cite{SW00} as a natural extension of the class of arithmetic Fuchsian groups which includes all triangle groups (and their subgroups of finite index). A cofinite Fuchsian group $\Gamma$ is called \emph{semi-arithmetic} if the subgroup $\Gamma^{(2)}$ generated by the squares of the elements of $\Gamma$ is contained in an arithmetic group $\Delta$ acting on a product $\Hyp^r$ of hyperbolic planes. We refer to Section~\ref{sec:prelim} for more details regarding the definition. By Takeuchi~\cite{Tak77}, almost all hyperbolic triangle groups are non-arithmetic but it is not hard to check that they are all semi-arithmetic. Other examples of strictly semi-arithmetic groups can be found in \cite{SW00, Kuch15, Cosac21, CD22}. Semi-arithmetic groups appear in the study of Teichm\"uller curves, billiards, translation surfaces and related topics. Lately these connections were investigated by \mbox{McMullen} in  \cite{McMullen20, McMullen21, McMullen22}. 

In this paper we are interested in finiteness properties of semi-arith\-me\-tic groups. Recall that by Borel's theorem there exist only finitely many arithmetic Fuchsian groups of bounded coarea. This is no longer the case for semi-arithmetic groups, as, for instance, there exist infinitely many hyperbolic triangle groups and all of them have bounded coarea. Therefore, in order to separate a finite subset of semi-arithmetic groups we need to consider other invariants. The second natural invariant is the \emph{arithmetic dimension} $r$ from the definition of a semi-arithmetic group. Indeed, Nugent and Voight showed in \cite{Nugent17} that there exist only finitely many triangle groups with bounded arithmetic dimension (see also Corollary~\ref{cor2} below for a different proof). Moreover, in Corollary~\ref{cor1} we show that an upper bound on coarea and arithmetic dimension imply finiteness of all semi-arithmetic groups that admit a modular embedding (by \cite{Kuch15} there exist Veech groups with modular embeddings that are not subgroups of hyperbolic triangle groups). This is about as far as one can get with these two invariants at hand -- it was shown in \cite{CD22} that for any fixed genus $g$ there exist semi-arithmetic surfaces defined over \emph{any} totally real number field of odd prime degree. In particular, there exist infinitely many semi-arithmetic groups of a fixed sufficiently large coarea with bounded arithmetic dimension. 

In order to have a finer control over the structure of semi-arithmetic groups we introduce a new invariant called \emph{stretch}. The name is inspired by Thurston's stretch of maps between hyperbolic surfaces from \cite{thurston1998minimal}, though our notion is different. Perhaps a more precise name would be \emph{arithmetic stretch} but we prefer to stay with a shorter version. Stretch is defined using Lipschitz geometry and the proof of its invariance with respect to commensurability is based on the notion of the Riemannian center of mass (see Definition~\ref{def stretch} and Proposition~\ref{prop:comm}). The Riemannian center of mass was introduced by Grove and Karcher in \cite{grove1973conjugate} and further developed in a series of papers by Karcher and collaborators in 1970s. It appeared to be extremely useful in differential geometry but we could not trace any prior applications of this construction in geometry of groups. We define stretch and discuss its basic properties in Section~\ref{sec:stretch}. 
The stretch of an arithmetic Fuchsian group is equal to one. The contracting property of holomorphic maps implies that the stretch of groups with modular embeddings is also equal to one. This way stretch can be understood as a quantitative measure of non-arithmeticity of a group. For us, it provides the missing parameter for the class of semi-arithmetic groups. 

We are now ready to state the main result of the paper:
\begin{thm}\label{thm1}
For any $L \geq 1$, $\mu>0$ and $r \ge 1$ there exist only finitely many conjugacy classes of semi-arithmetic Fuchsian groups with arithmetic dimension at most \(r\), stretch at most \(L\) and coarea at most \(\mu\).
\end{thm}

The proofs of Theorem~\ref{thm1} and related results together with some corollaries are presented in Section~\ref{sec:results}. The key new ingredient of the proof of the theorem is the \emph{arithmetic Margulis lemma}. This result was obtained in \cite{Fraczyk22} following the previous work of Breuillard in \cite{Breu11}. In order to prove Corollaries~\ref{cor1} and \ref{cor2} we use the contracting property of  modular embeddings which follows from the classical Schwarz--Pick lemma. 

As it was indicated above, all the three conditions in Theorem~\ref{thm1} are necessary. In particular, the results of \cite{CD22} imply that there exist infinite sets of semi-arithmetic groups with bounded arithmetic dimension and coarea whose stretch is unbounded. In Section~\ref{sec:examples} we modify the construction from \cite{CD22} in order to produce such groups defined over a \emph{fixed} totally real field $K$. As a corollary we obtain infinite families of non-commensurable Fuchsian groups contained in $\mathrm{SL}(2, K)$. The question about existence of such infinite sets was asked by A.~Rapinchuk and answered by Vinberg for $K = \mathbb{Q}$ \cite{vinberg2012}. Our result in Proposition~\ref{prop:examples} gives a construction of such sets of groups for the totally real fields $K$ of arbitrarily large degree. 

\medskip

\noindent
\textbf{Acknowledgements.} We thank Curtis McMullen for his interest in this work and helpful comments. 
The work of M.B. is partially supported by the FAPERJ research grant. G.C. is grateful for the grant 2022/10772-0, S\~ao Paulo Research Foundation (FAPESP). Part of this work was done while C.D. and G.T.P. were visiting IMPA supported by PRONEX.

\section{Preliminaries}\label{sec:prelim}

Let $\Gamma < \PSL(2,\R)$ be a finitely generated Fuchsian group. The \emph{trace field} of $\Gamma$ is the field generated by the traces of all elements of $\Gamma$ over the rational numbers $\Q(\mathrm{tr} \,\Gamma)$. This field, however, is not invariant under commensurability, which is the desirable notion of equivalence when dealing with subgroups of the integral points of some algebraic group (in particular, arithmetic and semi-arithmetic groups). We are then motivated to look at $\Gamma^{(2)} = \langle \gamma^2 \mid \gamma \in \Gamma \rangle$. This is a finite index normal subgroup of $\Gamma$ with the property of having the minimal trace field among all finite index subgroups of $\Gamma$. This implies, in particular, that $\Q(\tr\,\Gamma^{(2)})$ is an invariant of the commensurability class of $\Gamma$, it is called the \emph{invariant trace field} of $\Gamma$ and denoted by $k\Gamma$.

\begin{definition}
    Let $\Gamma$ be a Fuchsian group of finite covolume. We say that $\Gamma$ is \emph{semi-arithmetic} if $k\Gamma$ is a totally real number field and the traces of elements of $\Gamma$ are algebraic integers.
\end{definition}

This notion is closely related to \emph{arithmetic Fuchsian groups}. Indeed, arithmetic groups must also satisfy the following: every nontrivial Galois embedding $k\Gamma \hookrightarrow \R$ maps the set $\tr\,\Gamma^{(2)}$ to a bounded subset of $\R$ (cf. \cite{Takeuchi75} and \cite{maclachlan2003arithmetic}). In the semi-arithmetic case, these sets may be unbounded for, say, $r$ of the Galois embeddings. This number is called the \emph{arithmetic dimension} of $\Gamma$ (cf. \cite{Nugent17}). Note that, if $n=[k\Gamma:\Q]$, then $1\leq r \leq n$.  In particular, arithmetic Fuchsian groups have arithmetic dimension $1$.

Alternatively, arithmetic and semi-arithmetic Fuchsian groups may be characterised in terms of quaternion algebras, as we present next.

Let $A = \hilbert{a}{b}{k}$ be a quaternion algebra over the totally real number field $k$ of degree $[k:\Q] = n$. Denote by $\sigma_1 = \id, \sigma_2,\dots,\sigma_n$ the Galois embeddings of $k$ into $\R$. Each of these extends to an embedding $\rho_i$ of $A$ into $\hilbert{\sigma_i(a)}{\sigma_i(b)}{\R}$, which is a quaternion algebra over $\R$ and thus must be isomorphic either to the algebra $\rM(2,\R)$ of $2\times 2$ matrices with real coefficients, or to the algebra $\cH$ of Hamilton's quaternions. In the former case, we say that $A$ is \emph{unramified} over $\sigma_i$, or that $A$ \emph{splits} over $\sigma_i$. Otherwise, we say $A$ is \emph{ramified} over $\sigma_i$.

Assume $A$ splits over exactly $r\geq 1$ embeddings and suppose, without loss of generality, that they have been labeled so that $A$ splits over $\sigma_1,\dots,\sigma_r$. Then we obtain an isomorphism
\begin{align}\label{14723.1}
    A\otimes_\Q \R \cong  \rM(2,\R)^{r} \times \cH^{n-r},
\end{align}
by mapping $x\otimes a$ to $(a\rho_1(x),\dots,a\rho_n(x))$.

Note that the natural inclusion $A \hookrightarrow A\otimes_{\Q}\R$, followed by  the isomorphism \eqref{14723.1} and the projection onto the $i$th factor, gives back the embedding $\rho_i: A \hookrightarrow \rM(2,\R) \text{ or } \cH$. We observe that each $\rho_i$ preserves both norm and trace, in the sense that $\tr\,\rho_i(x) = \sigma_i(\tr(x))$ and $\det\rho_i(x) = \sigma_i(\rn(x))$, for all $x \in A$.  

In particular, if $A^1$ denotes the subgroup of elements of $A$ of norm $1$, then $\rho = (\rho_1,\dots,\rho_r)$ restricts to an embedding
\begin{align*}
    \rho: A^1 \hookrightarrow \SL(2,\R)^r.
\end{align*}
Note that, for each $x\in A^1$, $\rho(x) = (\rho_i(x),\dots,\rho_r(x))$ acts on $\bH^r$ componentwise via M\"obius transformations.

Let $\cO$ be an order in $A$. The well-known Borel--Harish-Chandra Theorem implies that $\rho$ maps $\cO^1$ onto a discrete subgroup of $\SL(2,\R)^r$ of finite covolume. In other words, $\rho(\cO^1)$ is a lattice in $\SL(2,\R)^r$.

\begin{definition}\label{031121.1}
    A subgroup $\Delta$ of $\PSL(2,\R)$ is said to be an \emph{arithmetic group acting on} $\bH^r$ if it is commensurable to some $\rP\rho_1(\cO^1)$ as above. In case $\Delta$ is a finite index subgroup of $\rP\rho_1(\cO^1)$, we say it is \emph{derived from a quaternion algebra}.
\end{definition}

This definition was given in \cite{SW00}. In particular, any group derived from a quaternion algebra as above, acts on $\bH^r$ in a natural way since it can be embedded into $\SL(2,\R)^r$ by the map $\rho \circ \rho_1^{-1}$. 

Let $\Gamma$ be a semi-arithmetic Fuchsian group. The subset of $\rM(2,\R)$ consisting of finite sums of the form $\sum a_i\gamma_i$ where each $\gamma_i$ is an element of $\Gamma^{(2)}$ and each $a_i$ is in $k\Gamma$ inherits a natural structure of an algebra over $k\Gamma$. In fact, it is a quaternion algebra over $k\Gamma$, denoted by $A\Gamma$, which is an invariant of the commensurability class of $\Gamma$ (see \cite[Chapter 3]{maclachlan2003arithmetic}). The subring $\cO\Gamma^{(2)}$, consisting of those finite sums whose coefficients are algebraic integers in $k\Gamma$, forms an order in $A\Gamma$. Then $(\cO\Gamma^{(2)})^1$ is an arithmetic group acting on $\bH^r$ (in fact, it is derived from a quaternion algebra), and $\Gamma$ is commensurable to the subgroup $\Gamma^{(2)}$ of $(\cO\Gamma^{(2)})^1$. The embedding $\rho$ maps $\Gamma^{(2)}$ to a subgroup of a lattice in $\SL(2,\R)^r$. This proves one direction of the following characterisation of semi-arithmetic Fuchsian groups. 

\begin{proposition}[{\cite[Proposition 1]{SW00}}]
    A Fuchsian group $\Gamma$ of finite covolume is semi-arithmetic if and only if it is commensurable to a subgroup of an arithmetic group acting on $\bH^r$.
\end{proposition}

We say that a semi-arithmetic Fuchsian group $\Gamma$ is derived from a quaternion algebra if $\Gamma < \Delta$, for some arithmetic group $\Delta $ acting on $\Hyp^r$ which is derived from a quaternion algebra.
It follows from the previous paragraph that \(\Gamma^{(2)}\) is always  derived from a quaternion algebra.

An important subclass of semi-arithmetic groups is defined using the notion of \emph{modular embedding}. It was first formally introduced in \cite{SW00}. We will use a slightly different definition which includes certain necessary adjustments (see \cite{MR1887671} and  \cite{McMullen23}).

\begin{definition}
    Let $\Gamma$ be a semi-arithmetic Fuchsian group so $\Gamma^{(2)}$ is derived from a quaternion algebra.  This provides a family of embeddings $f: \Gamma^{(2)} \to \SL(2,\R)^r$ (the restriction of $\rho\circ\rho_1^{-1}$ in the comment following Definition~\ref{031121.1}, for each choice of $\rho$).
    We say that $\Gamma$ \emph{admits modular embedding} 
    if $f$ can be chosen such that 
    there exists an $f$-equivariant holomorphic function $F: \bH \to \bH^r$, i.e., a function $F$ satisfying:
    \begin{align*}
        F(\gamma \cdot z) = f(\gamma)\cdot F(z),
    \end{align*}
    for all $z\in \bH$ and all $\gamma \in \Gamma^{(2)}$. 
\end{definition}

With this notion of modular embedding, Proposition~2(ii) in \cite{MR1887671} allows us to apply the argument of \cite{MR1075639} to prove the following important result.

\begin{theorem}\label{thm:trianglegroups}
All Fuchsian triangle groups admit modular embedding. 
\end{theorem}

The notion of arithmeticity is deeply related with commensurability of lattices. Recall that the \emph{commensurator} of a subgroup $\Gamma$ of a group $G$ is defined by 
$$ \mathrm{Comm}(\Gamma) = \{ g\in G \mid [\Gamma : \Gamma \cap g^{-1}\Gamma g] < \infty\}.$$ 
A fundamental result of Margulis \cite[Theorem~1, p.~2]{Margulis} states that a lattice $\Gamma$ in a semisimple Lie group is non-arithmetic if and only if its commensurator is discrete, and hence is a maximal lattice containing $\Gamma$. Notice that properly semi-arithmetic groups satisfy this property.  

We conclude this section with a technical lemma that will be used later on.  

\begin{lemma}\label{Lemma for conjugation}
 Let $\Gamma,\Gamma'$ be cofinite Fuchsian groups of the same signature such that $\Gamma$ is maximal and let $\phi:\Gamma \to \Gamma'$ be an isomorphism of groups. Consider a finite index normal subgroup $N \vartriangleleft \Gamma$ and define $N'=\phi(N)$. If the groups $N$ and $N'$ are conjugate, then $\Gamma$ and $\Gamma'$ are conjugate.
\end{lemma}
\begin{proof}
 Suppose that $N'=g N g^{-1}$ for some $g \in \PSL(2,\R)$. Since $\Gamma$ has the same coarea as $g^{-1} \Gamma' g$, the groups  $\Gamma$ and  $g^{-1} \Gamma' g$ are equal if and only if $g^{-1} \Gamma' g$ is contained in $\Gamma$. In order to prove this inclusion, it is sufficient to show that $g^{-1} \Gamma' g$ is contained in the \emph{normalizer} $\mathcal{N}_{\PSL(2,\R)}(N)=\{ \beta \in \PSL(2,\R) \mid \beta N \beta^{-1} = N \}$ of $N$. Indeed, $\mathcal{N}_{\PSL(2,\R)}(N)$ is a cofinite Fuchsian group (see \cite[Corollary 4.5.5]{Morris15}) containing $\Gamma$ and, since $\Gamma$ is maximal, we have $\Gamma=\mathcal{N}_{\PSL(2,\R)}(N)$. 
 
 Let $\gamma \in \Gamma$ and consider $\eta=g^{-1} \phi(\gamma) g \in g^{-1} \Gamma' g$. We have
 $$\eta N \eta^{-1}=(g^{-1} \phi(\gamma) g)(g^{-1} N' g)(g^{-1} \phi(\gamma^{-1}) g).$$
 
 Now we use that $N \vartriangleleft \Gamma$  and the definition of $N'$ in order to obtain 
 $$\eta N \eta^{-1}= g^{-1}\phi(\gamma N \gamma^{-1}) g= g^{-1}N' g= N.$$
\end{proof}

\section{Stretch of a semi-arithmetic group}\label{sec:stretch}

Let \(\Gamma\) be a semi-arithmetic Fuchsian group with arithmetic dimension \(r\) and commensurable to a subgroup of a group \(\Delta\) derived from a quaternion algebra. 

Let \(G=\mathrm{PSL}(2,\mathbb{R})^r\) be the orientation-preserving isometry group of \(\mathbb{H}^r\), and let \(\rho_0:\Delta \to G \) be a discrete representation such that \(\rho_0(\Delta)\) is an arithmetic subgroup. If \(\Gamma < \Delta\), then, up to replacing \(\Delta\) by a finite index subgroup, we can suppose that \(\Delta\) (hence \(\Gamma\)) is torsion-free. In this case,  \(S=\Gamma \backslash \mathbb{H}\) and \(M=\rho_0(\Delta)\backslash\mathbb{H}^r\) are smooth manifolds. The representation  \(\rho_0 \), when restricted to \(\Gamma\), defines a homomorphism \(\rho\) that can be seen as a homomorphism  from \(\pi_1(S)\) to \(\pi_1(M).\) Since \(M\) is an aspherical manifold, there exists a continuous map \(u:S \to M\) with \( \rho = u_*:\pi_1(S) \to \pi_1(M)\) (see \cite[Theorem 1.7.6]{Martelli16}). By the Whitney approximation theorem we can suppose that \(u\) is smooth. Moreover, we can suppose that such \(u\) is constant outside a compact set if \(S\) is noncompact. Hence, by lifting \(u\) we show that the set of \(\rho\)-equivariant smooth Lipschitz maps from \(\mathbb{H}\) to \(\mathbb{H}^r \) is nonempty. We are considering here the \emph{Kobayashi metric} on $\mathbb{H}^r$, i.e. the supremum of the hyperbolic metrics on each factor.   

If $\Gamma$, and therefore $\Delta$, is not torsion-free, we consider finite index subgroups $\Gamma' < \Gamma$ and $\Delta' < \Delta$ which are torsion free, and then let $\Gamma_0 = \Gamma' \cap \Delta' < \Gamma_1 = \Gamma \cap \Delta'< \Delta' $. By the previous argument, there exists a function $u_0:\mathbb{H} \to \mathbb{H}^r$ which is smooth, $\rho$-equivariant with respect to $\Gamma_0$ and $L$-Lipschitz for some  $L \in \mathbb{R}$. We now show that it is possible to construct a function $u:\mathbb{H} \to \mathbb{H}^r$ which is also smooth, $\rho$-equivariant with respect to $\Gamma_1$ and  $L$-Lipschitz.

First we notice that, for $j=1,\ldots,r $, if $\pi^j: G \to \PSL(2,\R)$ and $\pi_j: \Hyp^r \to \Hyp$,  are the  natural projections, we can write $u_0 = (u_0^1, \ldots, u_0^r)$ and $\rho = (\rho_1, \ldots,  \rho_r )$, where for each $j$, $u_0^j :=\pi_j \circ u_0$, $\rho_j = \pi^j \circ \rho$ and $u_0^j$ is $\rho_j $-equivariant with respect to $\Gamma_0$. Therefore, if we construct, for each $j$, a map $u_j: \Hyp \to \Hyp$, which is $\rho_j$-equivariant with respect to $\Gamma_1$ and $L$-Lipschitz, we get that $u= (u_1,\ldots, u_r): \Hyp \to \Hyp^r$ is the function that we are looking for.

Let $n = |\Gamma_1:\Gamma_0|$ and write $\Gamma_1 = g_1\Gamma_0 \cup\ldots\cup  g_n\Gamma_0 $. We fix $j \in \{1,\ldots,r\}$ and consider, for each \(i=1,\ldots,n\), the map 
$$F_i(z)=\rho_j(g_i)\cdot u_0^j(g_i^{-1}z).$$ 
The desired map $u_j$ will be the \emph{center of mass} of the \(F_i's\), which can be obtained by using a construction of Karcher (see Section 1 of \cite{karcher1977riemannian}). If we let \(B=\{1,\ldots,n\}\) with the usual counting measure, we can define for each \(z \in \mathbb{H}\) the mass distribution \(Z:B \to \mathbb{H}\) given by \(Z(i)=F_i(z)\). The map $$P_Z(x)=\frac{1}{2}\int_B d(x,Z(a))^2 da = \frac{1}{2n}\sum_{i=1}^n d(x,F_i(z))^2$$
is smooth by Theorem 1 in [op. cit.], and its gradient \(\nabla P_Z\) is given by \[\nabla P_Z(x)= \int_B \exp_x^{-1}(Z(a)) da=\frac{1}{n}\sum_{i=1}^n \exp_{x}^{-1}(F_i(z)).\]
Moreover, since \(\mathbb{H}\) has negative curvature, \(\nabla P_Z\) has a unique singularity which is, by definition, the center of mass \(u_j(z)\) of the set \(\{F_1(z),\ldots,F_n(z)\}\). By the Implicit Function Theorem, the function \(u_j(z)\) is smooth. Furthermore, by Theorem 1.5 in \cite{karcher1977riemannian}, we have \[\lvert \nabla (P_Z(x))\rvert \geq d(x, u_j(z)).\] 
This inequality implies the map $F$ defined above is $L$-Lipschitz. Indeed, we have
\begin{eqnarray*}
d(u_j(w), u_j(z)) &\leq & |\nabla (P_Z(u_j(w)))|\\
&=& \big\lvert \frac{1}{n}\displaystyle\sum_{i=1}^n\exp_{u_j(w)}^{-1}F_i(z)\big\rvert\\
&\leq & \frac{1}{n}\displaystyle\sum_{i=1}^n |\exp_{u_j(w)}^{-1}F_i(z) - \exp^{-1}_{u_j(w)}F_i(w)|\\
&\leq &\frac{1}{n}\displaystyle\sum_{i=1}^n d(F_i(z),F_i(w))\\
&\leq &Ld(z,w),
\end{eqnarray*}
where we used the following facts:
\begin{enumerate}
\item by the definition of $u_j(w)$ with \(W(i)=F_i(w)\), we have
 \[0=\nabla P_W(u_j(w))=\frac{1}{n}\displaystyle\sum_{i=1}^n \exp_{u_j(w)}^{-1}F_i(w);\]  
\item the function $\exp_p^{-1}$ is $1$-Lipschitz in the hyperbolic plane for any $p$;  and 
\item the functions $F_i$ are $L$-Lipschitz. 
\end{enumerate}

Therefore, there exists a function $u:\mathbb{H} \to \mathbb{H}^r$ which is smooth, $\rho$-equivariant with respect to $\Gamma_1$ and  $L$-Lipschitz.

\medskip

We are now ready to define the notion of \textit{stretch} of a semi-arithmetic Fuchsian group derived from a quaterion algebra with representation \(\rho_0\). We say that such a group has stretch at most \(L\) if there exists a \(\rho_0\)-equivariant \(K\)-Lipschitz map  with \(K\leq L\).

For the general case we recall that \(\Gamma^{(2)}\) is derived from a quaternion algebra. Therefore, the stretch of \(\Gamma\) can be defined as follows: 

\begin{definition}\label{def stretch}
A semi-arithmetic Fuchsian group \(\Gamma\) has \textit{stretch at most \(L\)} if  \(\Gamma^{(2)}\) has stretch at most \(L.\) This means that, for $r = \mathrm{a.dim}(\Gamma)$, there exist a discrete representation \(\rho_0:\Gamma^{(2)} \to \mathrm{PSL}(2,\mathbb{R})^r \) and a $\rho_0$-equivariant \(K\)-Lipschitz map  with \(K\leq L\).  The infimum of the set of such constants \(L\)  will be called stretch of the group $\Gamma$ and denoted by \(\delta(\Gamma)\). We also write $\delta(S) = \delta(\Gamma)$ for the stretch of the orbifold $S = \Gamma \backslash \mathbb{H}$. 
\end{definition}

\begin{proposition}\label{prop:comm}
    Let $L>0$ be a fixed number. The property of having stretch at most $L$ is invariant by commensurability.
\end{proposition}

\begin{proof}
Suppose \(\Gamma_1\) is commensurable to \(\Gamma_2\) and let $r$ be their arithmetic dimension. Then there exists a torsion-free group \(H\), which is a finite index subgroup of both \(\Gamma_1^{(2)}\) and \(\Gamma_2^{(2)}\). If $\Gamma_1$, and thus $\Gamma_1^{(2)}$ has stretch at most $L$, there is a $K$-Lipschitz map $u:\mathbb{H} \to \mathbb{H}^r$, for some $K\leq  L$, which is $\Gamma_1^{(2)}$-equivariant, and therefore $H$-equivariant. Hence $H$  has also stretch at most $L$.

By using the center of mass construction as above, since $H$ has finite index in $\Gamma_2^{(2)}$, we can construct a $K$-Lipschitz $\Gamma_2^{(2)}$-equivariant map, and therefore $\Gamma_2^{(2)}$, and thus $\Gamma_2$, also have stretch at most $L$.
\end{proof}

Stretch can be estimated in terms of the matrix coefficients of the group. In order to state this result we will need to recall some notations. Given a subgroup $\Gamma < \mathrm{PSL}(2, \mathbb{R})$, we denote by $\mathrm{Spec}(\Gamma)$ the set of the the biggest eigenvalues of the hyperbolic elements of $\Gamma$. If $\lambda$ is an algebraic integer, the \emph{house} $\ho{\lambda}$ is the maximal absolute value of its conjugates. 

\begin{proposition}
\label{prop:stretch}
     If \(\Gamma\) is a semi-arithmetic Fuchsian group with representation $\rho$ as above, which has stretch at most \(L\), then  
 \begin{equation}\label{eq-stretch-spec}
           \sup \left\{ \frac{ \log \ho{\lambda}}{\log \lambda} \mid \lambda \in \mathrm{Spec}(\Gamma^{(2)}) \right\} \leq L.
\end{equation}
\end{proposition}

\begin{proof}
    We can suppose that \(\Gamma\) is derived from a quaternion algebra. In this case, for any hyperbolic element \(\gamma \in \Gamma\) and \(z\) contained in the axis of \(\gamma,\) any \(L\)-Lipschitz \(\rho\)-equivariant map \(f: \Hyp \to \Hyp^r\) satisfies 
\[2\log(\lambda_i(\gamma))=\ell(\rho_i(\gamma)) \leq d(f_i(z),f_i(\gamma z)) \leq L \ell(\gamma) = 2L\log(\lambda_1(\gamma)),  \]    
for each \(i=1,\ldots,r\), 
where \(\lambda_i(\gamma)\) is the biggest eigenvalue of \(\rho_i(\gamma)\) if \(\rho_i(\gamma)\) is hyperbolic and 1 otherwise. 
Since \(\ho{\lambda(\gamma)}=\max\{\lambda_i(\gamma)\}\) and \(\lambda=\lambda_1(\gamma)\) we conclude that \(\log(\ho{\lambda(\gamma)}) \leq L \log(\lambda). \)
\end{proof}

By the definition, the stretch of an arithmetic Fuchsian group is equal to~$1$. The next proposition shows that the same holds true for all semi-arithmetic groups that admit modular embeddings.

\begin{proposition}\label{prop:stretch modular}
If a semi-arithmetic Fuchsian group $\Gamma$ admits a modular embedding, then its stretch is equal to $1$.
\end{proposition}

\begin{proof}
Since $\Gamma$ admits a modular embedding, there exists a holomorphic $\Gamma$-equivariant function $f:\mathbb{H} \to \mathbb{H}^r$. By the Schwarz--Pick lemma, this function is $1$-Lipschitz and thus $\delta(\Gamma) \leq 1$. On the other hand, the previous proposition implies that $\delta(\Gamma) \ge 1$.
\end{proof}

In Section~\ref{sec:examples} we will construct examples of semi-arithmetic groups $\Gamma$ with arbitrarily large stretch. Let us note in passing that we did not find the exact value of stretch for these groups but only a lower bound. Indeed, we do not have a single example of a group for which the stretch is bigger than $1$ and computed explicitly. We leave it as an open problem.

\begin{problem}
Find the exact value of $\delta(\Gamma)$ for some semi-arithmetic group $\Gamma$ with $\delta(\Gamma) > 1$.
\end{problem}

Another interesting question arises in connection with Proposition~\ref{prop:stretch}: \textit{Is it possible to give an upper bound for $\delta(\Gamma)$ in terms of the spectrum? Or, perhaps, even define the stretch by the left hand side of \eqref{eq-stretch-spec}?} 
For our main applications we use both the lower spectral bound and the upper geometrical bound. This allows us to turn stretch into an effective invariant and to prove new results. Our previous investigations were focused around the spectral properties of semi-arithmetic groups but without the geometric counter-part their applications were limited.

We conclude this section with a brief discussion of some other notions related to stretch of embeddings. We have already mentioned that the notion of stretch was motivated by the term used by Thurston in the context of maps between surfaces \cite{thurston1998minimal}. In this paper we generalize it for maps from semi-arithmetic surfaces to associated higher dimensional arithmetic locally symmetric spaces. In our definition of stretch and the proof of its commensurability invariance we used Karcher's center of mass. 
Another related geometric construction is the well known Douady--Earle barycentric extension for homeomorphisms of the circle to the disk \cite{DE86}. Although the ideas are similar (in both cases the goal is to find an average of points in order to define an equivariant map), the methods are quite different. In particular, in \cite{DE86}, the authors define a map that associates a barycentre in the Poincar\'e disc to a measure on a circle, while we are interested in higher dimensional spaces of non-constant sectional curvature. The results of Douady and Earle were greatly generalized by Besson, Courtois and Gallot (cf. \cite{BCG96}) and it may be possible to apply their approach to our maps. One important feature of Karcher's method is that it allows us to construct a new map with the same Lipschitz constant. This property is essential for our applications and it is not clear if it can be achieved via the Besson--Courtois--Gallot method. We leave a detailed study of this question for future research. 


\section{Main results}\label{sec:results}

For any given integers \(g, n\geq 1 \) with \(g+n \geq 2,\) let $\mathcal{M}_{g,n}$ be the moduli space of surfaces with genus $g$ and $n$ cusps. In addition, for any integer \(r \geq 1\) and positive number \(L \geq 1\), define 
\[\mathcal{SA}_{g,n}(r,L)=\{ S \in \mathcal{SA}_{g,n} \mid  \mathrm{a.dim}(S) \leq r \mbox{ and }  \delta(S) \leq L\},\]
where \(\mathcal{SA}_{g,n}=\{ S \in \mathcal{M}_{g,n} \mid S \mbox{ is semi-arithmetic}\}\) and \(\mathrm{a.dim}(S)\) denotes the arithmetic dimension of \(S.\)

\begin{theorem}\label{corollary from Margulis}
  Let \(\Gamma\) be a torsion-free semi-arithmetic Fuchsian group with invariant trace field \(k\), arithmetic dimension \( \leq r\),  stretch \(\leq L\) and coarea \( \leq \mu\). Then
   \[[k:\mathbb{Q}] \leq C \log(\mu)+c\]
for some constants \(C,c>0\) which depend only on \(r\) and \(L.\) 
\end{theorem}
\begin{proof}
    If \(\Gamma\) has parabolic elements, then \(A\otimes_\mathbb{Q}\mathbb{R}\) has no compact factors and hence \([k:\mathbb{Q}]=r\). Therefore, the theorem holds with the constants \(C=1\) and \(c=r-\log(\frac{\pi}{3})\), using that  $\frac{\pi}{3}$ is the minimal coarea of a non-cocompact Fuchsian group.

    If \(\Gamma\) is cocompact of coarea \(\mu(\Gamma)\), by Yamada's theorem (see the work of Yamada \cite{yamada1982marden} or the proof given in Theorem 5.3 of \cite{fanoni2015maximum} for a more precise result) there exist a point \(p\in \mathbb{H}\) and a constant \(R>0\) such that the set 
    \[D(p,R)=\{g \in \Gamma \mid d(gp,p) \leq R\}\]
    contains two primitive non-conjugate hyperbolic elements. For all  \(g \neq 1\), we have \(g B(p,R) \cap B(p,R)= \emptyset\), hence 
    \[R \leq \cosh^{-1}\left(\frac{\mu(\Gamma)}{2\pi}+1\right).\]
    Moreover, the subgroup of \(\Gamma \) generated by \(D(p,R)\) contains a nonabelian free group.

  Since \(\Gamma^{(2)}\) is derived from a quaternion algebra and torsion-free, this group is embedded in a group \(\Delta\) derived from a quaternion algebra acting on \(\mathbb{H}^r\) and there must exist a \(L\)-Lipschitz  \(\Gamma^{(2)}\)-equivariant map \(u:\mathbb{H} \to \mathbb{H}^r.\) Thus, if \(g \in D(p,R)\), we get that \(g^2 \in \Gamma^{(2)}\) and 
  \[d_r(g^2 \cdot u(p),u(p))=d_r(u(g^2 \cdot p),u(p)) \leq 2LR.\]

 By the arithmetic Margulis lemma in \cite[Section~3]{Fraczyk22}, there exists $\epsilon_r >0$, which depends only on $r$, such that the subgroup generated by
\begin{align}\label{261122.1}
    \{\lambda \in \Delta \mid d_r(\lambda \cdot u(p),u(p)) < \epsilon_r [k:\mathbb{Q}] \}
\end{align}
is virtually nilpotent. Since the subgroup generated by \(\{g^2 \mid g \in D(p,R)\}\) contains a nonabelian free group we obtain
\[\epsilon_r [k:\mathbb{Q}] \leq 2LR. \]
Therefore,
\[[k:\mathbb{Q}] \leq C\log(\mu)+c\]
for some constants \(C,c>0\) which depend only on \(r\) and \(L.\)
\end{proof}

\begin{proposition}\label{lower bound for systole}
Let  \(S\) be a semi-arithmetic Riemann surface of co\-area \(\mu\), arithmetic dimension \(r\) and stretch at most \(L\), then \(\sys(S) \geq s\) for some  positive constant \(s=s(\mu,r,L)\).

\end{proposition}    
\begin{proof}
Let \(S=\Gamma \backslash \Hyp\), where \(\Gamma \) is a semi-arithmetic group with invariant trace field \(k\) and invariant quaternion algebra \(A\). Consider an arbitrary hyperbolic element $\gamma \in \Gamma$ with biggest eigenvalue \(\lambda(\gamma)>1\). Note that \(\gamma^2\) satisfies the quadratic equation \(x^2-tx+1\) in \(A\), with \(t=\lambda(\gamma)^2+\lambda(\gamma)^{-2} \in R_k\), where $R_k$ is the ring of integers of $k$.

If we denote by $P(x)$ the minimal polynomial of $\lambda(\gamma)^2$, then $P(x)$ is an integral monic polynomial of degree $D=2[k:\mathbb{Q}]$. We can assume that \(D \geq 4\), because \(k=\mathbb{Q}\) reduces to the arithmetic case. By Theorem~\ref{corollary from Margulis}, $D \leq 2C\log(\mu)+2c$, where \(C\) and \(c\) are positive constants depending only on \(r\) and \(L\).

Let $\Mah(P)$ be the \emph{Mahler measure} of $P(x)$, which is given by 
\[ \Mah(P) = \prod_{i=1}^D \mathrm{max}(1, |\theta_i|),\]
where $\theta_1$,\ldots, $\theta_D$ are the roots of $P(x)$.

Since \(A\) splits at \(r\) real places and \(S\) has stretch at most \(L\), we obtain by Proposition~\ref{prop:stretch} that
\[\Mah(P)  \leq \lambda(\gamma)^{2rL}. \]

On the other hand, by Dobrowolski's bound for the Mahler measure \cite{dobrowolski1979question} we have $$ \log (\Mah(P)) \ge U \left(\frac{\log (\log(D)}{\log(D)} \right)^3 $$ for some universal constant $U>0.$ 

Hence, $$ \ell(\gamma)=\log(\lambda(\gamma)^2) \geq \frac{\log(\Mah(P))}{rL} \ge \frac{U}{rL} \left(\frac{\log (\log(D))}{\log(D)} \right)^3.$$

Since $\dfrac{\log(\log(x))}{\log(x)}$ is decreasing for $x > 1,$ and \(D \leq 2C\log(\mu)+2c\), there exists a constant \(s=s(\mu,r,L) > 0\) which depends only on \(\mu,r\) and \(L\) such that \(\ell(\gamma) \geq s\) for any hyperbolic element \(\gamma \in \Gamma\) .
\end{proof}    

\begin{theorem}\label{Torsion-free}
The set \(\mathcal{SA}_{g,n}(r,L)\) is finite, for all positive integers \(r,g\) and  \(n\) and every \(L \geq 1\).
\end{theorem}
\begin{proof}
    Given \(S=\Gamma \backslash \mathbb{H} \in \mathcal{SA}_{g,n}(r,L)\), by Theorem~\ref{corollary from Margulis} the invariant trace field of \(\Gamma\) has  bounded degree \(d\).

  By Proposition~\ref{lower bound for systole}, there exists \(\delta > 0\) depending only on \((g,n)\), \(r\) and \(L\) such that \(\sys(S) \geq \delta.\) Hence, by Theorem 1 in \cite{bers1972remark}, \(S\) is contained in a compact subset \(\mathcal{K}\) of \(\mathcal{M}_{g,n}\). 
  
  Let \(\ell_1,\ldots,\ell_s\) be a set of length functions defined on a fundamental domain of the modular space \(\mathcal{M}_{g,n}\) in its universal covering which determine uniquely a surface \( S \). If \(\ell_i(P) \leq \ell\) for any \(i=1,\ldots,s\) and for any \(P \in \mathcal{K}\), then the map which associates \(S\) to the parameters \(\left(\exp(\frac{\ell_i(S)}{2})\right)_i\) is an injective map from \(\mathcal{SA}_{g,n}(r,L)\) to \(R_{2d}(\ell L)^s\) which is finite, where \(R_n(X)\) is the set of algebraic integers of degree at most \(n\) and house at most \(X.\) 
\end{proof}

We are now ready to prove Theorem~\ref{thm1}, which we restate for convenience.

\begin{theorem}[Main Theorem] \label{thm:main}
For any $L \geq 1$, $\mu>0$ and $r \ge 1$ there exist finitely many conjugacy classes of semi-arithmetic Fuchsian groups with arithmetic dimension at most \(r\), stretch at most \(L\) and coarea at most \(\mu\).
\end{theorem}

\begin{proof}
We can suppose that the Fuchsian groups are maximal. For any conjugacy class of a semi-arithmetic Fuchsian group satisfying the hypothesis of the theorem we can associate its signature \(\sigma\). Then we can split the proof in two parts: \(i)\) proving that there are finitely many possible signatures, and \(ii)\) showing that for a given signature there are finitely many conjugacy classes of maximal semi-arithmetic Fuchsian groups with bounded arithmetic dimension and stretch.

For part \( (i)\) we recall that the area \(\mu\) of a surface \( S = \Gamma \backslash \mathbb{H}\) is expressed in terms of the signature of the Fuchsian group \( \Gamma \). Therefore, the genus, the number of cusps, and the number of conjugacy classes of elliptic elements of a group \(\Gamma\) of signature \(\sigma\) are bounded above by a function depending only on \(\mu.\) It remains to show that the maximal order of elliptic elements is also bounded.

By Theorem 1.2 of \cite{edmonds1982torsion}, there exists a torsion-free subgroup \( \Gamma_1 <\Gamma\) of index at most \(2t^m\), where \(t\) (resp. \(m\)) is the maximal order (resp. the number of conjugacy classes) of elliptic elements of \(\Gamma\). Since \(\Gamma\) and \(\Gamma_1\) share the same invariant trace field \(k\), they have the same arithmetic dimension and both have stretch at most \(L\). As \(\Gamma_1\) has coarea at most \(2t^m \mu\), by Theorem~\ref{corollary from Margulis} we have
\[[k:\mathbb{Q}] \leq Cm\log(t)+C\log(2\mu)+c \leq C'\log(t)\]
for some constant \(C'\) which depends only on \(r\), \(L\) and \(\mu\).

On the other side, since a quadratic extension of \(k\) contains \(2\cos(\frac{2\pi}{t}),\) we have that \([k:\mathbb{Q}] \geq \frac{1}{2}\phi(t),\) where \(\phi\) denotes Euler's totient function. The inequality \[\frac{1}{2}\phi(t) \leq C'\log(t)\]
can hold only for finitely many values of \(t\), hence we conclude the proof of \( (i)\).

For part \( (ii)\) we take the torsion-free group \( \Gamma_1 \) as above and consider the embedding of marked conjugacy classes of maximal Fuchsian groups of signature \(\sigma\) satisfying the hypothesis of theorem into the set of marked conjugacy classes of Fuchsian groups of signature of the group \(\Gamma_1\) (cf. \cite{greenberg1963maximal}).  By Theorem~\ref{Torsion-free}, up to conjugation there are only finitely many possibilities for \(\Gamma_1\), and since \(\Gamma\) is maximal we can apply Lemma~\ref{Lemma for conjugation} in order to guarantee that up to conjugation there are finitely many possibilities for \(\Gamma.\) 
\end{proof}

By Proposition~\ref{prop:stretch modular}, we have the following consequence of Theorem~\ref{thm:main}:

\begin{corollary}\label{cor1}
    For any $\mu>0$ and $r \ge 1$ there exist finitely many conjugacy classes of semi-arithmetic Fuchsian groups which admit modular embeddings, with arithmetic dimension at most \(r\) and coarea at most~\(\mu\).
\end{corollary}

In \cite{Nugent17}, Nugent and Voight proved that the number of triangle groups with bounded arithmetic dimension is finite. Theorem~\ref{thm:main} provides another proof of this fact, due to Theorem \ref{thm:trianglegroups} and Corollary~\ref{cor1} we have:

\begin{corollary}\label{cor2}
    For any $r \ge 1$ there exist finitely many conjugacy classes of triangle groups with arithmetic dimension at most \(r\).
\end{corollary}

\section{Infinite sequences of semi-arithmetic groups with the same quaternion algebra
}
\label{sec:examples}

It is natural to ask if the conditions in Theorem~\ref{thm:main} are all necessary. The existence of infinitely many non-commensurable triangle groups shows the necessity to bound arithmetic dimension. There are 
 also infinitely many non-commensurable arithmetic groups with arbitrary large coarea. We now show that there exists an infinite family of non-commensurable semi-arithmetic Fuchsian groups with the same coarea and arithmetic dimension, but with unbounded stretch. 

In what follows, let \(R_K\) be the ring of integers of the real quadratic number field \(K = \mathbb{Q}(\sqrt{3})\)  and let \(\epsilon \in R_K^*\) be a fixed nontrival unit, from what we get \(\mathbb{Q}(\epsilon)=K\).

Let us fix \(x=\sinh^{-1}(\frac{\epsilon}{2})\). There exists a unique (up to isometry) trirectangle $\mathcal{Q}$ with acute angle equal to \(\frac{\pi}{3}\) and one of the opposite sides to this angle of length \(x\). Let \(y\) (resp. \(z\)) denote the length of the other opposite side to the acute angle (resp. the length of the diagonal joining two right angled vertices of the trirectangle). By the hyperbolic trigonometry of the trirectangle we obtain 
\begin{align}
\sinh(x)\sinh(y)= &\frac{1}{2} \label{eq.1}, \\
\cosh(x)\cosh(y)= & \cosh(z).\ \label{eq.2}
\end{align}


\vspace{-2cm}
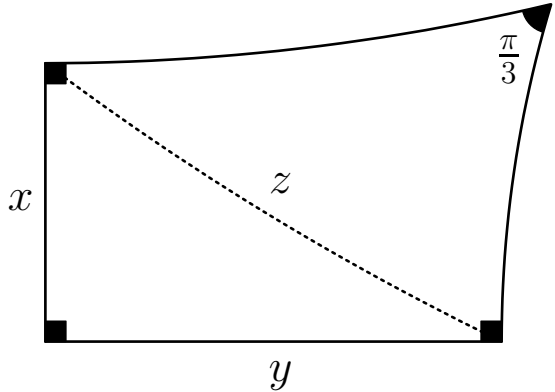
\begin{figure}[!ht]
    \centering
    \begin{tikzpicture}[scale = 15, line cap=round,line join=round,>=triangle 45,x=1cm,y=1cm]

        \clip(-0.1,-0.05) rectangle (0.5,0.5);
        
        \draw[line width=0.5pt,fill=black,fill opacity=0.10000000149011612] (0.4,0.018055031292605973) -- (0.38194496870739403,0.018055031292605977) -- (0.38194496870739403,0) -- (0.4,0) -- cycle; 
        \draw[line width=0.5pt,fill=black,fill opacity=0.10000000149011612] (0,0.22794496870739403) -- (0.018055031292605977,0.22794496870739403) -- (0.018055031292605973,0.246) -- (0,0.246) -- cycle; 
        \draw[line width=0.5pt,fill=black,fill opacity=0.10000000149011612] (0.018055031292605973,0) -- (0.01805503129260598,0.01805503129260597) -- (0,0.018055031292605973) -- (0,0) -- cycle; 
        \draw [shift={(0.4433153434209247,0.2981100831367445)},line width=0.5pt,fill=black,fill opacity=0.10000000149011612] (0,0) -- (193.4161859031799:0.025533670123074) arc (193.4161859031799:253.4853804916076:0.025533670123074) -- cycle;

        \draw [line width=1pt] (0,0)-- (0.4,0);
        \draw [line width=1pt] (0,0)-- (0,0.246);
        \draw [shift={(1.45,-0.0003629205329044032)},line width=1pt]  plot[domain=2.8533581684542586:3.141247015000791,variable=\t]({1*1.050000062719671*cos(\t r)+0*1.050000062719671*sin(\t r)},{0*1.050000062719671*cos(\t r)+1*1.050000062719671*sin(\t r)});
        \draw [shift={(0.00026327611383750144,2.1555203252032515)},line width=1pt]  plot[domain=4.712251104862395:4.946545597454828,variable=\t]({1*1.9095203433529173*cos(\t r)+0*1.9095203433529173*sin(\t r)},{0*1.9095203433529173*cos(\t r)+1*1.9095203433529173*sin(\t r)});
        \draw [shift={(1.450893183196884,2.1569726556046898)},line width=1pt,dotted]  plot[domain=4.062998359169128:4.259027672242962,variable=\t]({1*2.399355646734353*cos(\t r)+0*2.399355646734353*sin(\t r)},{0*2.399355646734353*cos(\t r)+1*2.399355646734353*sin(\t r)});

        \draw[color=black] (-0.028,0.12537589864964804) node {\larger[2] $x$};
        \draw[color=black] (0.2,-0.03) node {\larger[2] $y$};
        \draw[color=black] (0.2,0.14) node {\larger[2] $z$};
        \draw[color=black] (0.4,0.25) node {\larger[2] $\frac{\pi}{3}$};

    \end{tikzpicture}
    
    \caption{Trirectangle}\label{191023.1}
\end{figure}

Consider the index two subgroup \(\Gamma\) of orientation-preserving isometries of the reflection group whose fundamental domain is $\mathcal{Q}$. This Fuchsian group is generated by some elements \(X,Y,Z,W\) such that \(X,Y,Z\) are elliptic of order two, \(W\) has order three and \(XYZW=1\). We note that this group is generated by the  elements \(A=XY\) and \(B=ZX\). Indeed,
\begin{align*}
 [A,B]=(XY)(ZX)(YX)(XZ)=(XYZ)(XYZ)=W^{-2}=W.   
\end{align*}
Moreover,
\begin{align*}
     X=(YZW)^{-1}=W^2(ZY)=W^2(ZX)(XY)=W^2AB. 
\end{align*}

Therefore $\Gamma$ is a non-elementary Fuchsian group generated by \(A\) and \(B\). By Lemma~3.5.7 in \cite{maclachlan2003arithmetic} the invariant trace field of \(\Gamma\) is 
\begin{align*}
    k_\Gamma=\mathbb{Q}(\tr(A)^2,\tr(B)^2,\tr(A)\tr(B)\tr(AB)).
\end{align*}
It follows from Theorem 3.6.2 [op. cit.] that the invariant quaternion algebra \(A_\Gamma\) can be described as
\begin{align*}
    A_\Gamma= \left( \frac{\tr(A)^2(\tr(A)^2-4),\tr(A)^2\tr(B)^2(\tr([A,B])-2)}{k_\Gamma}   \right).
\end{align*}

We now compute $ k_\Gamma$ and $ A_\Gamma$ explicitly. First we note that \(A\) and \(B\) are products of order two elliptic elements, hence they are hyperbolic, and their displacements are given by twice the distance between the fixed points of the respective half-turns (see Theorem 7.38.2 in \cite{Beardon12}). Similarly, \(BA\) is hyperbolic, with displacement given by twice the diagonal of the trirectantle. 

We can suppose that \(A\) has displacement \(2x\). In this case \(B\) has displacement \(2z\). If \(\lambda\) (resp., \(\mu\) and \(\eta\)) denotes the biggest positive eigenvalue of \(A\) (resp., \(B\) and \(BA\)), then \(\lambda=e^x\), ~\(\mu=e^z\) and \(\eta=e^y\). Since \(\sinh(x)=\frac{\epsilon}{2}\) and by \eqref{eq.1}, \(\sinh(y)=\epsilon^{-1}\), we have
\[\lambda-\lambda^{-1}=2\sinh(x)=\epsilon \mbox{ and } \eta-\eta^{-1}=2\epsilon^{-1}.\]
Hence
\begin{align} 
 \lambda &  = \frac{\epsilon+\sqrt{4+\epsilon^2}}{2}  \label{eq.3} ,\\
\lambda^{-1} & = \frac{-\epsilon+\sqrt{4+\epsilon^2}}{2} .& \label{eq.4}
\end{align}

Analogously,
\begin{align}
 \eta &  =\epsilon^{-1}+\sqrt{1+\epsilon^{-2}} \label{eq.5} ,\\
\eta^{-1} & =-\epsilon^{-1}+\sqrt{1+\epsilon^{-2}}.
\end{align}

By definition, \(\tr(A)=\lambda+\lambda^{-1}\), hence
\begin{align*}
\tr(A)^2(\tr(A)^2-4) & =(  \lambda^2+\lambda^{-2}+2)(\lambda^2+\lambda^{-2}-2) \\
  & = \lambda^4+\lambda^{-4} - 2   \\ &= (\lambda^2-\lambda^{-2})^2 \\
 & =\epsilon^2(4+\epsilon^2). 
\end{align*}
where the last equality follows from \eqref{eq.3} and \eqref{eq.4}.

By \eqref{eq.2}, we obtain
$$\tr(B)=\mu+\mu^{-1} = \frac{1}{2}(\lambda+\lambda^{-1})(\eta+\eta^{-1}) = \sqrt{4+\epsilon^2}\sqrt{1+\epsilon^{-2}}.$$

In order to calculate \(k_\Gamma\), note that 
\begin{align*}
 \tr(A)\tr(B)\tr(AB) & = (\lambda+\lambda^{-1})(\mu+\mu^{-1})(\eta+\eta^{-1})   \\
 & = \frac{1}{2}(\lambda+\lambda^{-1})^2(\eta+\eta^{-1})^2 \\
 & = \frac{1}{2}(4+\epsilon^2)(1+\epsilon^{-2}).
\end{align*}

Therefore,  
$k_\Gamma=\mathbb{Q}(4+\epsilon^2,(4+\epsilon^2)(1+\epsilon^{-2}))=\mathbb{Q}(\epsilon)=K.$

Since \([A,B]\) has order three, we get \(\tr([A,B])=\pm 1\) which implies that \(\tr([A,B])-2 = -1 \cdot \tau^2\) , for \(\tau \in K\). Hence, by the general property that \( \left(\displaystyle \frac{au^2,bv^2}{k} \right) \cong \left(\displaystyle \frac{a,b}{k} \right)\) for any \(u,v \in k\), we have
\[A_\Gamma=\left(\frac{4+\epsilon^2,-(1+\epsilon^{-2})}{k_\Gamma} \right) \cong \left(\frac{4+\epsilon^2,-(1+\epsilon^{2})}{k_\Gamma} \right).\]

We use Theorem $4.3.1$ in \cite{maclachlan2003arithmetic} to conclude that $A_\Gamma \cong \mathrm{M}(2,K)$, because the equation $(4+\epsilon^2)x^2 - (1+\epsilon^{2})y^2 = 1$ has  $x=y = \frac{1}{\sqrt{3}} $ as a  solution in $K$.

We claim that $\Gamma$ is semi-arithmetic. Indeed, $k_{\Gamma}$ is a totally real number field. From the computation above, we obtain that the traces of $A$, $B$ and $AB$ are algebraic integers. Thus by Lemma~3.5.2 in \cite{maclachlan2003arithmetic}, the traces of the elements of $\Gamma$ are algebraic integers.

For every nontrivial unit of the form $\epsilon = \epsilon_0^n$, for $n\in \mathbb{N}$, where $\epsilon_0$ is the fundamental unit of the field, we can reproduce the same construction. This gives a semi-arithmetic Fuchsian group $\Gamma_n$ with invariant trace field $K$ and invariant quaternion algebra $\mathrm{M}(2,K)$.
Notice that none of the groups groups $\Gamma_n$ is arithmetic, as their arithmetic dimension is equal to the degree of $K$, which is 2. Each \(\Gamma_n\) has coarea \(\frac{\pi}{3}\)  and arithmetic dimension \(2\). Moreover, the stretch \(\delta(\Gamma_n)\) goes to infinity as $n$ grows: by using the same notation as above, if we fix a unit \(\epsilon=\epsilon_0^n\), the corresponding hyperbolic element \((BA)^2 \in \Gamma_n^{(2)}\) has  \(\tau_n=\eta^2=1+2\epsilon^{-2}+2\epsilon^{-1}\sqrt{1+\epsilon^{-2}}\) as its biggest eigenvalue. This is an algebraic integer of degree four whose conjugates are \(\tau_n,\tau_n^{-1},\omega_n\) and \(\omega_n^{-1}\), where \(\omega_n=\ho{\tau_n}=1+2\epsilon^{2}+2\epsilon\sqrt{1+\epsilon^{2}}\). Thus, by Proposition~\ref{prop:stretch}, we have the estimate: 
$$\delta(\Gamma_n) \geq \frac{\log(\ho{\tau_n})}{\log(\tau_n)}.$$
Since \(\ho{\tau_n} \to \infty\) and \(\tau_n \to 1\), we get \(\lim \delta(\Gamma_n) \to \infty.\) As the stretch is invariant by conjugation, there must exist infinitely many of the groups $\Gamma_n$ which are pairwise non-conjugate.

The same construction works for any totally real number field $K$ which contains $\mathbb{Q}(\sqrt{3})$.

It is worth mentioning that in \cite{vinberg2012}, Vinberg gave an answer to a question by A.~Rapinchuk about the existence of an infinite family of non-commensurable (up to conjugacy) non-arithmetic Fuchsian groups with the same quaternion algebra. Another construction providing an infinite family of non-commensurable Fuchsian groups with the same algebra was given later by Norfleet  \cite{Norfleet2015}. The groups constructed by Vinberg and Norfleet sit in $\mathrm{SL}(2, \mathbb{Q})$, they have arithmetic dimension one and are quasi-arithmetic. The question about existence of infinite families of lattices in the
group of quaternions of norm $1$ of quaternion algebras over other fields remained open. Our construction allows us to obtain many examples of this kind:

\begin{proposition}\label{prop:examples}
    Given a totally real field extension $K/\mathbb{Q}(\sqrt{3})$ of degree $d\ge 1$, there exist infinitely many non-commensurable semi-arithmetic cocompact Fuchsian groups  contained in $\mathrm{PSL}(2,K)$. 
\end{proposition}

\begin{proof}
Given a field $K$ and the algebra $A = \mathrm{M}(2,K)$, we consider the  
groups $\Gamma_n$ constructed above. Pairwise non-commensurability of infinitely many of these semi-arithmetic groups follows immediately from Proposition~\ref{prop:comm}. The non-commensurability can also be checked by a more direct argument. To this end we note that  each group $\Gamma_n$ is a maximal lattice. Indeed, since $\Gamma_n$ is a non-arithmetic two generator Fuchsian group, we can apply Proposition~2.1 from \cite{MR91} to conclude that its commensurator, which is the maximal element in its commensurability class, must be either the group $\Gamma_n$ itself or a triangle group. However, by the complete description of triangle groups that contain two-generator subgroups given in [op. cit.], no non-arithmetic group of signature $(0;2,2,2,3;0)$ (which is the case of $\Gamma_n$) is contained in a triangle group. Therefore, the groups $\Gamma_n$ are maximal Fuchsian groups. They are pairwise non-conjugate, hence by maximality and non-arithmeticity they must be non-commensurable.
\end{proof}


\bibliographystyle{alpha}
\bibliography{Bibliography}

\end{document}